\documentclass[10pt,a4paper,reqno]{amsart}

\usepackage{latexsym,floatflt}
\usepackage{epsfig}
\usepackage{rotating}
\usepackage{amssymb}
\usepackage[T1]{fontenc}
\usepackage{afterpage}
\usepackage{color}

\usepackage{amssymb,amsfonts,amsmath,stmaryrd}

\catcode`\@=11
\def\section{\@startsection{section}{1}%
 \z@{.7\linespacing\@plus\linespacing}{.5\linespacing}%
 {\normalfont\bfseries\scshape\centering}}
\def\subsection{\@startsection{subsection}{2}%
  \z@{.5\linespacing\@plus\linespacing}{.5\linespacing}%
  {\normalfont\bfseries\scshape}}

\def\subsubsection{\@startsection{subsubsection}{3}%
 \z@{.5\linespacing\@plus\linespacing}{-.5em}
  {\normalfont\bfseries\itshape}}
\catcode`\@=12

\def\qed{$\hfill{\vrule height 3pt width 5pt depth 2pt}$}

\newfont{\bbold}{msbm9 scaled \magstep1}
\newfont{\bbolds}{msbm7 scaled \magstep1}

\newcommand{\ns}{\mbox{\bbold N}}

\newcommand{\zs}{\mbox{\bbold Z}}
\newcommand{\zss}{\mbox{\bbolds Z}}
\newcommand{\qs}{\mbox{\bbold Q}}
\newcommand{\rs}{\mbox{\bbold R}}

\newcommand{\bx}{\bar x}

\newcommand{\GK}{\mathbb{K}}

\newcommand{\C}{\mathcal C}
\newcommand{\E}{\mathcal E}
\newcommand{\M}{\mathcal M}

\newcommand{\Ref}[1]{(\ref{#1})}
\newcommand{\beq}{\begin{equation}}
\newcommand{\eeq}{\end{equation}}
\newcommand{\gf}{generating function}
\newcommand{\gfs}{generating functions}

\newcommand{\eps}{\epsilon}

\def\NN{\mathbb{N}}

\def\xx{\bar x}
\def\yy{\bar y}\def\zz{\bar z}

\def\ZZ{\mathbb{Z}}

\def\ct{\mathop{\mathrm{CT}}}

\def\pt{\mathop{\mathrm{PT}}}
\def\nt{\mathop{\mathrm{NT}}}

\def\sy{\mathfrak{S}}

\def\emm#1,{{\em #1}}

 \newtheorem{Theorem}{Theorem}
 \newtheorem{Proposition}[Theorem]{Proposition}

\newtheorem{Lemma}[Theorem]{Lemma}

\title[On partitions avoiding 3-crossings]
{On partitions avoiding 3-crossings}

\author{Mireille Bousquet-M\'elou}
\address{CNRS, LaBRI, Universit\'e Bordeaux 1, 351 cours de la Lib\'eration,
  33405 Talence Cedex, France}
\email{mireille.bousquet@labri.fr}

\author{Guoce Xin}
\address{Department
of Mathematics, Brandeis University,
Waltham MA 02454-9110, USA}
\email{guoce.xin@gmail.com}
\date
{January 20, 2006}

\begin{document}
\maketitle

\begin{flushright}
{To Xavier Viennot, on the occasion of his 60th birthday}
\end{flushright}

\begin{abstract}
A partition on $[n]$ has a crossing if there exists
$i_1<i_2<j_1<j_2$ such that $i_1$ and $j_1$ are in the same block,
$i_2$ and $j_2$ are in the same block, but $i_1$ and $i_2$ are not
in the same block. Recently, Chen et al. refined this classical
notion by introducing $k$\emm-crossings,, for any integer $k$. In
this new terminology, a classical crossing is a 2-crossing. The
number of partitions of $[n]$ avoiding $2$-crossings is well-known
to be the $n$th Catalan number $C_n={{2n}\choose n}/(n+1)$. This
raises the question of counting $k$-noncrossing partitions for $k\ge
3$. We prove that the sequence counting 
$3$-noncrossing
partitions is P-recursive, that is, satisfies a linear recurrence
relation with polynomial coefficients. We give explicitly such a
recursion. However, we conjecture that $k$-noncrossing partitions
are not P-recursive, for $k\ge 4$.

We obtain similar results for partitions avoiding \emm enhanced,
$3$-crossings. 
\end{abstract}

\section{Introduction}

A partition of $[n]:=
\{1,2, \ldots , n\}$ is a
collection of nonempty and mutually disjoint subsets of $[n]$,
called \emph{blocks}, whose union is $[n]$. The number of partitions of
$[n]$ is the \emph{Bell
  number} $B_n$.  A well-known
refinement of $B_n$ is given by the \emph{Stirling number} (of the second
kind) $S(n,k)$. It counts partitions of $[n]$ having
exactly $k$ blocks.

Recently another refinement of the Bell numbers by \emm crossings, and
\emm nestings, has attracted some
interest~\cite{crossings-nestings,kratti-diagrams}. This 
refinement is based on
the standard representation of a partition $P$ of $[n]$ by
a graph, whose vertex set $[n]$ is identified with the points $i\equiv (i,0)$
on the plane, for $1\le i \le n$, and whose edge set  consists of arcs
connecting 
the elements that occur \emm consecutively, in the same block (when each
block is totally ordered). For example, the standard representation of
$1478-236-5$ is given by the following graph.

\medskip
\begin{center}
\begin{picture}(0,0)%
\includegraphics{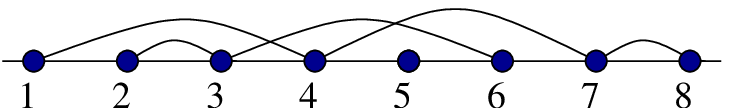}%
\end{picture}%
\setlength{\unitlength}{3947sp}%
\begingroup\makeatletter\ifx\SetFigFont\undefined%
\gdef\SetFigFont#1#2#3#4#5{%
  \reset@font\fontsize{#1}{#2pt}%
  \fontfamily{#3}\fontseries{#4}\fontshape{#5}%
  \selectfont}%
\fi\endgroup%
\begin{picture}(3474,487)(739,-1336)
\end{picture}%

\end{center}
\smallskip

Then  crossings and nestings have a  natural definition. A
$k$-\emph{crossing} of $P$ is a collection of $k$ edges
$(i_1,j_1)$, $(i_2,j_2)$, \dots, $(i_k,j_k)$
 such that $i_1<i_2<\cdots < i_k <j_1<j_2<\cdots <j_k$. This
means a subgraph of $P$ as drawn as follows.

\medskip
\begin{center}
\begin{picture}(0,0)%
\includegraphics{kcrossing.pstex}%
\end{picture}%
\setlength{\unitlength}{3947sp}%
\begingroup\makeatletter\ifx\SetFigFont\undefined%
\gdef\SetFigFont#1#2#3#4#5{%
  \reset@font\fontsize{#1}{#2pt}%
  \fontfamily{#3}\fontseries{#4}\fontshape{#5}%
  \selectfont}%
\fi\endgroup%
\begin{picture}(4252,595)(739,-2294)
\put(1201,-2236){\makebox(0,0)[lb]{\smash{{\SetFigFont{12}{14.4}{\familydefault}{\mddefault}{\updefault}{\color[rgb]{0,0,0}$i_2$}%
}}}}
\put(826,-2236){\makebox(0,0)[lb]{\smash{{\SetFigFont{12}{14.4}{\familydefault}{\mddefault}{\updefault}{\color[rgb]{0,0,0}$i_1$}%
}}}}
\put(2176,-2236){\makebox(0,0)[lb]{\smash{{\SetFigFont{12}{14.4}{\familydefault}{\mddefault}{\updefault}{\color[rgb]{0,0,0}$i_k$}%
}}}}
\put(3151,-2236){\makebox(0,0)[lb]{\smash{{\SetFigFont{12}{14.4}{\familydefault}{\mddefault}{\updefault}{\color[rgb]{0,0,0}$j_1$}%
}}}}
\put(3676,-2236){\makebox(0,0)[lb]{\smash{{\SetFigFont{12}{14.4}{\familydefault}{\mddefault}{\updefault}{\color[rgb]{0,0,0}$j_2$}%
}}}}
\put(4576,-2236){\makebox(0,0)[lb]{\smash{{\SetFigFont{12}{14.4}{\familydefault}{\mddefault}{\updefault}{\color[rgb]{0,0,0}$j_k$}%
}}}}
\end{picture}%

\end{center}
\smallskip

\noindent
A different notion of $k$-crossing is obtained
by  representing each block  by a complete
graph
\cite[p.~85]{klazar}. A $k$-\emph{nesting} of $P$ is
a collection of $k$ edges
$(i_1,j_1)$, $(i_2,j_2)$, \dots, $(i_k,j_k)$ such that $i_1<i_2<\cdots
< i_k <j_k<j_{k-1}<\cdots <j_1$, as represented below.

\medskip
\begin{center}
\begin{picture}(0,0)%
\includegraphics{knesting.pstex}%
\end{picture}%
\setlength{\unitlength}{3947sp}%
\begingroup\makeatletter\ifx\SetFigFont\undefined%
\gdef\SetFigFont#1#2#3#4#5{%
  \reset@font\fontsize{#1}{#2pt}%
  \fontfamily{#3}\fontseries{#4}\fontshape{#5}%
  \selectfont}%
\fi\endgroup%
\begin{picture}(4252,795)(739,-2294)
\put(1201,-2236){\makebox(0,0)[lb]{\smash{{\SetFigFont{12}{14.4}{\familydefault}{\mddefault}{\updefault}{\color[rgb]{0,0,0}$i_2$}%
}}}}
\put(826,-2236){\makebox(0,0)[lb]{\smash{{\SetFigFont{12}{14.4}{\familydefault}{\mddefault}{\updefault}{\color[rgb]{0,0,0}$i_1$}%
}}}}
\put(2176,-2236){\makebox(0,0)[lb]{\smash{{\SetFigFont{12}{14.4}{\familydefault}{\mddefault}{\updefault}{\color[rgb]{0,0,0}$i_k$}%
}}}}
\put(3151,-2236){\makebox(0,0)[lb]{\smash{{\SetFigFont{12}{14.4}{\familydefault}{\mddefault}{\updefault}{\color[rgb]{0,0,0}$j_k$}%
}}}}
\put(4126,-2236){\makebox(0,0)[lb]{\smash{{\SetFigFont{12}{14.4}{\familydefault}{\mddefault}{\updefault}{\color[rgb]{0,0,0}$j_2$}%
}}}}
\put(4576,-2236){\makebox(0,0)[lb]{\smash{{\SetFigFont{12}{14.4}{\familydefault}{\mddefault}{\updefault}{\color[rgb]{0,0,0}$j_1$}%
}}}}
\end{picture}%

\end{center}
\smallskip

\noindent A partition is $k$-\emph{noncrossing} if it has no
$k$-crossings, and $k$-\emph{nonnesting} if it has no
$k$-nestings.

A  variation of $k$-crossings (nestings), called
\emph{enhanced} $k$-crossings (nestings), was also studied
in~\cite{crossings-nestings}. One first adds a loop to every isolated
point in the standard representation of partitions. Then by
allowing $i_k=j_1$ in a $k$-crossing, we get an enhanced
$k$-crossing; in particular, a partition avoiding enhanced
$2$-crossings has parts of size 1 and 2 only.  Similarly, by allowing
$i_k=j_k$ in the definition of a $k$-nesting, we get an enhanced
$k$-nesting.

Chen et al. gave in~\cite{crossings-nestings} a bijection between
partitions of $[n]$ and certain ``vacillating'' tableaux of length
$2n$.  Through this bijection, a partition is $k$-noncrossing if
and only if the corresponding tableau has
height less than $k$,
and $k$-nonnesting if the tableau has
width less than $k$. A
simple symmetry on tableaux then entails that $k$-noncrossing
partitions of $[n]$ are equinumerous with $k$-nonnesting
partitions of $[n]$, for all values of $k$ and $n$. A second bijection
 relates partitions to certain  ``hesitating'' tableaux, in such a
 way the size of the largest \emm enhanced, crossing (nesting) becomes the
 height (width) of the tableau. This implies that partitions of $[n]$
 avoiding enhanced $k$-crossings are equinumerous with partitions of  $[n]$
 avoiding enhanced $k$-nestings.

\medskip

The number $C_2(n)$ of $2$-noncrossing partitions of $[n]$ (usually called
\emm noncrossing partitions,~\cite{kreweras-noncroisees,rodica}) is well-known
to be  the  Catalan number
$C_n=\frac{1}{n+1}\binom{2n}{n}$. For  $k>2$, the number of $k$-noncrossing
partitions of $[n]$ is not known, to the extent of our
knowledge. However, for any $k$, the number of  $k$-noncrossing
\emm matchings, of $[n]$ (that is, partitions in which all blocks have
size 2) is known to form a P-recursive sequence, that is,  to satisfy
a linear recurrence relation with polynomial
coefficients~\cite{grabiner,crossings-nestings}.
In this paper, we  enumerate
$3$-noncrossing partitions of $[n]$ (equivalently, $3$-nonnesting
partitions of $[n]$).  We obtain for $C_3(n)$ a (not so simple) closed
form expression (Proposition~\ref{closed-form}) and a linear
recurrence relation with polynomial coefficients.

\begin{Proposition}\label{main-thm}
The number $C_3(n)$  of  $3$-noncrossing partitions is given by $
C_3\left( 0 \right)=C_3\left(
1 \right) =1,$ and for $n\ge 0$,
\begin{equation}\label{e-precursive}
9n \left( n+3 \right) C_3 \left( n \right) -2\, \left(
5{n}^{2}+32 n+42 \right) C_3 \left( n+1 \right) + \left( n+7
\right)  \left( n+6
 \right) C_3\left( n+2 \right)
=0  .
\end{equation}
Equivalently, the associated \gf \ $\C(t)=\sum_{n\ge 0} C_3(n)t^n$
satisfies
\begin{equation}\label{e-dfinite}
 t^2 (1-9t)(1-t) {\frac {d^{2}}{d{t}^{2}}}\C(t)
+ 2t\left( 5- 27{t}+ 18{t}^{2} \right) {\frac {d }{dt}}\C(t)
+ 10 \left(2 -3t \right) \C(t)-20
=0.
\end{equation}
Finally,  as
$n$ goes to infinity,
$$
C_3(n) \sim  \frac {3^9\cdot 5}{2^5}\frac {\sqrt 3}\pi \,\frac {9^n}{n^7}.
$$
\end{Proposition}
\noindent
The first few values of the sequence $C_3(n)$, for $n \ge 0$, are
$$
1, 1, 2, 5, 15, 52, 202, 859, 3930, 19095, 97566,\dots .
$$
A standard study~\cite{ince,wimp} of the above differential equation,
which can be done automatically  using the {\sc Maple} package {\sc DEtools},
suggests that $C_3(n) \sim \kappa \, {9^n}/{n^7}$
for some positive constant $\kappa$. However, one needs the explicit
expression of $C_3(n)$ given in Section~\ref{section-lagrange}
to prove this statement and find the value of $\kappa$.
The above asymptotic behaviour is confirmed experimentally
by the computation of the first values of
$C_3(n)$ (a few thousand  values can be computed  rapidly using the
package {\sc Gfun} of {\sc Maple}~\cite{Gfun}). For instance, when
$n=50000$,
then $n^7C_3(n)/9^n\simeq 1694.9$, while $\kappa\simeq 1695.6$.

 As discussed in the last section of the paper, the above result might
 remain isolated, as there is no (numerical) evidence that the \gf\ of
 4-noncrossing partitions should be P-recursive.

\medskip
We obtain a similar result for partitions avoiding \emm enhanced, $3$-crossings
(or enhanced 3-nestings). The number of partitions of $[n]$ avoiding
enhanced 2-crossings is easily seen to be the $n$th \emm Motzkin
number,~\cite[Exercise 6.38]{stanley-vol2}.

\begin{Proposition}\label{main-thm-H}
The number $E_3(n)$  of partitions of $[n]$ having no enhanced
$3$-noncrossing is given by $ E_3\left( 0 \right)=E_3\left( 1 \right)
=1,$ and for
$n\ge 0$,
$$8
\left( n+3 \right)  \left( n+1 \right) E_3 \left( n \right) +
 \left( 7{n}^{2}+53n+88 \right) E_3 \left( n+1 \right) - \left( n+8
 \right)  \left( n+7 \right) E_3\left( n+2 \right)=0
.
$$
Equivalently, the associated \gf \ $\E(t)=\sum_{n\ge 0} E_3(n)t^n$
satisfies
\begin{equation*}\label{e-dfinite-H}
{t}^{2} \left(1+ t \right)  \left(1- 8t \right) {\frac
{d^{2}}{d{t} ^{2}}}\E ( t )
+2t  \left( 6-23t-20{t}^{2} \right) {\frac {d}{dt}}\E ( t )
+6 \left(5-7t -4{t}^{2}\right) \E ( t ) -30=0.
\end{equation*}
Finally,  as
$n$ goes to infinity,
$$
E_3(n) \sim  \frac {2^{16}\cdot 5}{3^3}\frac {\sqrt 3}\pi \,\frac {8^n}{n^7}.
$$
\end{Proposition}
\noindent The first few values of the sequence $E_3(n)$, for $n
\ge 0$, are
$$
1, 1, 2, 5, 15, 51, 191, 772, 3320, 15032, 71084, \dots .
$$
Observe that $C_3(5)$ and $E_3(5)$ differ by 1: this difference comes
from the partition $135-24$ which has an enhanced 3-crossing but no
3-crossing (equivalently, from the partition $15-24-3$ which has an
enhanced 3-nesting but no 3-nesting).
Again, the study of the differential equation
suggests that $
E_3(n)\sim \kappa \frac{8^n}{n^7},
$
for some positive constant $\kappa$, but we  need the explicit
expression~\eqref{E-expr} of $E_3(n)$
to prove this statement and
find the value of $\kappa$. Numerically, we have found that for
$n=50000$, $n^7E_3(n)/8^n\simeq 6687.3$, while $\kappa \simeq 6691.1$.

\smallskip
The starting point of our proof of Proposition~\ref{main-thm} is
the above mentioned  bijection between partitions avoiding
$k+1$-crossings and vacillating tableaux of height at most $k$.
As described in~\cite{crossings-nestings}, these tableaux
can be easily encoded by certain $k$-dimensional lattice walks.
Let $D$ be a subset of $\ZZ^k$. A $D$-\emph{vacillating} lattice
walk of length $n$
is a sequence of lattice points $p_0,p_1,\dots ,p_{n}$ in $D$
such that for all $i$,
\begin{itemize}
\item [$(i)$] $p_{2i+1}=p_{2i}$ or $p_{2i+1}=p_{2i}-e_j$ for some
 coordinate vector $e_j=(0, \ldots, 0, 1, 0, \ldots ,0)$,
\item [$(ii)$] $p_{2i} = p_{2i-1}$ or
$p_{2i}=p_{2i-1} + e_j$ for some
$e_j$.
\end{itemize}
 We will be interested in two different domains $D$ of $\ZZ^k$: the domain
 $Q_k=\NN^k$ of points with
non-negative integer coordinates and the Weyl chamber (with
a slight change of coordinates) $W_k=\{\, (a_1,a_2,\dots
,a_{k})\in \ZZ^{k}: a_1> a_2> \cdots > a_{k}\ge 0 \,\}$. Vacillating
 walks in $W_k$ are related to $k+1$-noncrossing partitions as follows.

\begin{Theorem}[{\bf Chen et al.~\cite{crossings-nestings}}]
\label{t-chenetal}
Let $C_k(n)$ denote the number of $k$-noncrossing partitions of
$[n]$. Then $C_{k+1}(n)$ equals the number of
$W_k$-vacillating lattice walks of length $2n$ starting and ending at $(k-1,k-2,\dots,0)$.
\end{Theorem}

The proof of Proposition~\ref{main-thm} goes as follows: using the
reflection principle, we first reduce the enumeration of
vacillating walks in the Weyl chamber $W_k$ to that of vacillating
walks in the non-negative domain $Q_k$. This reduction is valid
for any $k$. We then focus on the case $k=2$. We write a
functional equation satisfied by  a 3-variable series that counts
vacillating walks in $Q_2$. This equation is based on a simple
recursive construction of the walks. It is  solved using a
2-dimensional version of the so-called \emm kernel method.,
This gives
the \gf \ $\C(t)=\sum C_3(n) t^n$ as the constant term in a certain algebraic
series. We then  use the Lagrange inversion formula to find an
explicit expression of $C_3(n)$, and apply
the
\emm creative telescoping, of~\cite{AB} to obtain the recurrence
relation. We finally derive from the expression of $C_3(n)$ the
asymptotic behaviour of these numbers.
The proof of
Proposition~\ref{main-thm-H}, given in Section~\ref{sec-hesitant}, is similar.

\section{Partitions with no $3$-crossing}
\subsection{The reflection principle}
Let $\delta=(k-1,k-2,\dots,0)$ and let $\lambda$ and $\mu$ be two
lattice points in $W_k$. Denote by $w_k(\lambda,\mu,n)$ (resp.
$q_k(\lambda,\mu,n)$) the number of $W_k$-vacillating (resp.
$Q_k$-vacillating) lattice walks of length $n$ starting at
$\lambda$ and ending at $\mu$. Thus Theorem \ref{t-chenetal}
states that $C_{k+1}(n)$ equals $w_k(\delta,\delta,2n)$.
 The reflection principle, in the vein of~\cite{gessel-zeil,zeil-andre},
 gives  the  following: 
\begin{Proposition}
\label{propo-reflection} For any starting and ending points $\lambda$ and
  $\mu$ in $W_k$, the number of $W_k$-vacillating walks going from
  $\lambda$ to $\mu$ can be expressed in terms of the number of
  $Q_k$-vacillating walks as follows:
\begin{align}\label{e-reflection}
w_k(\lambda,\mu,n)=\sum_{\pi\in \sy_k} (-1)^\pi
q_k(\lambda,\pi(\mu),n),
\end{align}
 where $(-1)^\pi$ is the sign of $\pi$ and
$\pi(\mu_1,\mu_2,\dots,\mu_k)=(\mu_{\pi(1)},\mu_{\pi(2)},\dots,
\mu_{\pi(k)})$.
\end{Proposition}

\noindent{\bf Proof.} Consider the set of hyperplanes $\mathcal{H}=\{\,
x_i=x_j : 1\le i<j \le k \, \}$. The reflection of the point $(a_1,
\ldots , a_k)$ with respect to the hyperplane $x_i=x_j$ is simply
obtained by exchanging the coordinates $a_i$ and $a_j$. In particular,
the set of (positive) steps taken by vacillating walks, being
$\{e_1, \ldots , e_k\}$, is invariant under
such reflections. The same holds for the negative steps, and of course
for the ``stay'' step, 0. This implies that reflecting a
$Q_k$-vacillating walk with respect to $x_i=x_j$ gives another
$Q_k$-vacillating walk. Note that this is not true
when reflecting with respect to $x_i=0$ (since $e_i$ is
transformed into $-e_i$).

Define a total ordering on  $\mathcal{H}$. Take a
$Q_k$-vacillating walk $w$ of length $n$ going from $\lambda$ to
$\pi(\mu)$, and assume it touches at least one hyperplane in
$\mathcal{H}$. Let $m$ be the first time  it touches a hyperplane.
Let $x_i=x_j$ be the smallest hyperplane it touches at time $m$.
Reflect all steps of $w$ after time $m$ across $x_i=x_j$; the
resulting walk $w'$ is a $Q_k$-vacillating walk going from
$\lambda$ to $(i,j)(\pi(\mu))$, where $(i,j)$ denotes the
transposition that exchanges $i$ and $j$. Moreover, $w'$ also
(first) touches $\mathcal{H}$ at time $m$, and the smallest
hyperplane it touches at this time is $x_i=x_j$.

The above transformation is  a sign-reversing involution on the set of
$Q_k$-vacillating paths that go from $\lambda$ to $\pi(\mu)$, for some
permutation $\pi$, and hit one of the hyperplanes of
$\mathcal{H}$. In the right-hand side of~\eqref{e-reflection}, the
contributions of these walks
  cancel out. One is left with the walks that  stay within the Weyl
  chamber, and this
happens only when $\pi$ is the identity. The proposition follows. \qed

\medskip
This proposition, combined with Theorem~\ref{t-chenetal}, gives
the number of $(k+1)$-noncrossing partitions of $n$ as a linear
combination of the numbers $q_k(\delta,\pi(\delta),2n)$.
Hence, in order to count $(k+1)$-noncrossing partitions, it
suffices to find a formula for $q_k(\delta,\mu,2n)$, for certain
ending points $\mu$. This is what we do below  for $k=2$.

\subsection{A functional equation}
\label{sec-functional}
Let us specialize Theorem~\ref{t-chenetal} and
Proposition~\ref{propo-reflection} to $k=2$. This gives
\begin{eqnarray}
C_3(n)&=& w_2((1,0),(1,0),2n)\nonumber \\
&=& q_2((1,0),(1,0),2n)-q_2((1,0),(0,1),2n).\label{C3-q}
\end{eqnarray}
From now on,  a lattice walk always means a  $Q_2$\emm-vacillating
lattice walk starting at, $(1,0)$, unless specified otherwise.

Let $a_{i,j}(n):=q_2((1,0),(i,j),n)$ be the number of lattice walks
of length $n$ ending at $(i,j)$. Let
$$F_e(x,y;t)=\sum_{i,j,n\ge 0} a_{i,j}(2n)x^iy^jt^{2n}$$
and
$$F_o(x,y;t)=\sum_{i,j,n\ge 0} a_{i,j}(2n+1)x^iy^jt^{2n+1}$$
be respectively the
 generating functions of lattice walks of even and odd length.
These series are power series in $t$ with
coefficients  in $\qs[x,y]$.
We will often work in a slightly larger ring, namely the ring
$\qs[x,1/x,y,1/y][[t]]$ of power series in $t$ whose coefficients are
Laurent polynomials in $x$ and $y$.

Now we find functional equations for $F_e(x,y;t)$
and $F_o(x,y;t)$.
By appending to an even length walk a west step
$(-1,0)$, or a south step $(0,-1)$, or a stay step $(0,0)$, we
obtain either an odd length walk (in $Q_2$), or an even length walk
ending on the
$x$-axis followed by a south step, or an even length walk ending
on the $y$-axis followed by a west step. This correspondence is easily seen
to be a bijection, and gives the following functional equation:
\begin{align}\label{e-Fe}
F_e(x,y;t)(1+\xx+\yy)t=F_o(x,y;t)+\yy tH_e(x;t)+\xx t V_e(y;t),
\end{align}
where $\xx= 1/x$, $\yy=1/y$, and $H_e(x;t)$ (resp.~$V_e(y;t)$) is the
generating function of even
lattice walks ending on the $x$-axis (resp.~on the $y$-axis).

Similarly, by adding to an odd length walk an east step $(1,0)$, or
a north step $(0,1)$, or a stay step, we obtain an even length
walk of positive length. The above correspondence is a bijection and gives
another functional equation:
\begin{align}\label{e-Fo}
F_o(x,y;t)(1+x+y)t=F_e(x,y;t)-x.
\end{align}
Solving equations \eqref{e-Fe} and \eqref{e-Fo} for $F_e(x,y;t)$
and $F_o(x,y;t)$ gives
\begin{align*}
F_e(x,y;t)&=\frac{x-t^2\yy(1+x+y)H_e(x;t)-t^2\xx(1+x+y)V_e(y;t)}
{1-t^2(1+x+y)(1+\xx+\yy)},\\
F_o(x,y;t)&= t\;\frac{x(1+\xx+ \yy)- \xx V_e(y;t)-\yy H_e(x;t)}
{1-t^2(1+x+y)(1+\xx+\yy)}.
\end{align*}
Since we are mostly interested in even length walks ending at $(1,0)$ and at
$(0,1)$, it suffices to determine the series $H_e(x;t)$ and $V_e(y;t)$,
and hence to solve one of the above functional equations. We
choose the one for $F_o(x,y;t)$, which is simpler.
\begin{Proposition}
The \gf\ $F_o(x,y;t)$ of $Q_2$-vacillating lattice walks of odd length
starting from $(1,0)$  is related to the \gfs\ of even lattice walks
of the same type ending on the $x$- or $y$-axis by
\begin{align}\label{e-Fkernel}
K(x,y;t)F(x,y;t) =xy+x^2y+x^2-x H(x;t)- y V(y;t),
\end{align}
where $F_o(x,y;t)=tF(x,y;t^2)$, $V(y;t)=V_e(y;t^2)$,
$H(x;t)=H_e(x;t^2)$, and $K(x,y;t)$ is the \emm kernel, of the equation:
\beq \label{kernel}
K(x,y;t)=xy-t(1+x+y)(x +y +xy).
\eeq
\end{Proposition}
From now on, we will very often
omit the variable $t$ in our notation. For instance, we will write
$H(x)$ instead of $H(x;t)$. Observe that~\eqref{e-Fkernel} defines
$F(x,y)$ uniquely as a series in $t$ with coefficients in
$\qs[x,y]$. Indeed, setting $y=0$ shows that $H(x)=x+t(1+x)F(x,0)$
while setting $x=0$ gives
$V(y)=t(1+y)F(0,y)$.

\subsection{The kernel method}
\label{sec-kernel}
We are going to apply to~\Ref{e-Fkernel} the {\em
obstinate kernel method\/} that has already been used
in~\cite{bousquet-motifs,bousquet-kreweras} to solve similar equations.
The classical kernel method consists in coupling the
variables $x$ and $y$ so as to cancel the  kernel $K(x,y)$.
This  gives some ``missing'' information about
the series $V(y)$ and $H(x)$ (see for
instance~\cite{bousquet-petkovsek-recurrences,hexacephale}). In its
obstinate version, the kernel
method is combined with a procedure that constructs and exploits
several  (related) couplings $(x,y)$. This procedure is essentially
borrowed from~\cite{fayolle-livre}, where similar functional equations
occur in a probabilistic context.

Let us start with the standard kernel method. First fix $x$, and
consider the kernel as a quadratic polynomial in $y$. Only one of
its roots, denoted $Y_0$ below, is a power series in $t$:
\begin{align} Y_0&= \frac
{1-(\xx+3+{x})t- \sqrt {\left(1-(1+x+ \xx)t\right)^2-4t}
}{ 2\left( 1+\xx \right) t} 
\nonumber \\
&=(1+x)t+\frac{(1+x)(1+3x+x^2)}{x}t^2+\cdots \nonumber
\end{align}
The coefficients of this series are Laurent polynomials in $x$, as is
easily seen from the equation
\beq
\label{eq-Y}
Y_0=t(1+x+Y_0)\left(1+(1+\xx)Y_0\right).
\eeq
Setting $y=Y_0$ in \eqref{e-Fkernel} gives a  valid identity between
series of $\qs[x,\xx][[t]]$, namely
$$
xH(x)+ Y_0V(Y_0)=xY_0+x^2Y_0+x^2.
$$
The second root of the kernel is $Y_1=x/Y_0= O(t^{-1})$, so that  the
expression $F(x,Y_1)$  is not well-defined.

\begin{figure}[hbt]
\begin{center}
\input{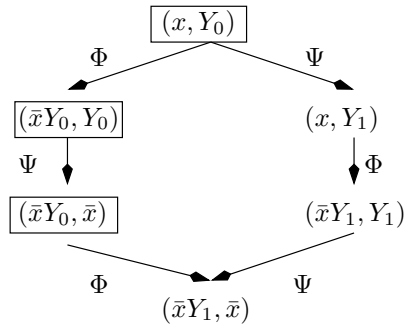}
\end{center}
\caption{The orbit of  $(x,Y_0)$ under the action of  $\Phi$ and $\Psi$.}
\label{diagram}
\end{figure}

 Now let $(X,Y)\not = (0,0)$ be a
pair of Laurent series in $t$ with coefficients in a field $\GK$
such that $K(X,Y)=0$. Recall that $K$ is quadratic in $x$ and $y$.
In particular, the equation $K(x,Y)=0$ admits a second solution
$X'$.
 Define $\Phi(X,Y)= (X',Y)$.
Similarly, define $\Psi(X,Y)= (X,Y')$, where $Y'$ is
the second solution  of $K(X,y)=0$. Note that $\Phi$ and $\Psi$ are
involutions. Moreover, with the kernel given by~\Ref{kernel}, one has
$Y'=X/Y$ and $X'=Y/X$. Let us examine the action of $\Phi$ and $\Psi$ on the
pair $(x,Y_0)$: we obtain an orbit of cardinality $6$
(Figure~\ref{diagram}).

The 6  pairs of power series given in Figure~\ref{diagram} cancel the
kernel, and we have framed
the ones that can be legally substituted for $(x,y)$ in the main functional
equation~\Ref{e-Fkernel}. Denoting $Y\equiv Y_0$, we thus obtain
 {\em three\/} equations relating the unknown series $H(x)$ and $V(y)$:

\begin{eqnarray}
xH(x)+ YV(Y)&=&xY+x^2Y+x^2,
\label{e-eq1}\\
\xx Y H(\xx Y)+ YV(Y)&=&\xx Y^2+\xx^{2}Y^3+\xx^{2} Y^2,
\label{e-eq2}\\
 \xx Y H(\xx Y)+ \xx V(\xx)&=&\xx^{2}Y+\xx^{3}Y^2+\xx^{2}Y^2.
\label{e-eq3}
\end{eqnarray}

\subsection{Positive and negative parts of power series}
A simple linear combination of the above three equations (namely,
\eqref{e-eq1}$-$\eqref{e-eq2}$+$\eqref{e-eq3}) allows us to eliminate
the terms $V(Y)$ and $ H(\xx Y)$. We are left with:
\begin{align*}
xH(x)+\xx V(\xx)=x^2+ (\xx^{2}+x+x^2)Y+(\xx^{3}-\xx) Y^2-\xx^{2} Y^3.
\end{align*}
Since $xH(x)$ contains only positive powers of $x$ and $\xx
V(\xx)$ contains only negative powers of $x$, we have characterized
the series $H(x)$ and $V(y)$.
\begin{Proposition}
\label{propo-walks}
The series $H(x)$ and $V(y)$  counting  $Q_2$-vacillating
walks of even
length starting at $(1,0)$ and ending on the $x$-axis and on the
$y$-axis satisfy 
\begin{align*}
xH(x)= \pt_x \left( x^2+ (\xx^{2}+x+x^2)Y+(\xx^{3}-\xx) Y^2-\xx^{2} Y^3\right)
,\\
\xx V(\xx)= \nt_x \left( x^2+ (\xx^{2}+x+x^2)Y+(\xx^{3}-\xx) Y^2-\xx^{2}
Y^3\right),
\end{align*}
where the operator $\pt_x$ (resp.  $\nt_x$) extracts positive
(resp. negative) powers of  $x$ in  series of $\qs[x,\xx][[t]]$.
\end{Proposition}
\noindent One may then go back to~(\ref{e-Fe}--\ref{e-Fo}) to
obtain expressions of the series $F_e(x,y;t)$ and $F_o(x,y;t)$.
However, our main concern in this note is the number $C_3(n)$ of
$3$-noncrossing partitions of $[n]$. Going back to \eqref{C3-q},
we find that $C_3(n)$ is determined by the following three
equations:
\begin{align*}
C_3(n)&  =q_2((1,0),(1,0),2n)-q_2((1,0),(0,1),2n),\\
q_2((1,0),(1,0),2n)&=[xt^n]\ H(x)=[x^2t^n] \ xH(x), \\
q_2((1,0),(0,1),2n)&=[yt^n]\ V(y)=[\xx^{2}t^n]\ \xx V(\xx).
\end{align*}
Using Proposition~\ref{propo-walks}, we obtain the
generating function $\C(t)$ of $3$-noncrossing partitions:
\begin{align}\label{e-Gtfirst}
\C(t)=\ct_x \; (\xx^{2}-x^2)\left(x^2+ (\xx^{2}+x+x^2)Y+(\xx^{3}-\xx)
Y^2-\xx^{2} Y^3\right),
\end{align}
where the operator $\ct_x$ extracts the constant term in $x$ of series
in $\qs[x,\xx][[t]]$.
Observe that $Y(x)=x Y(\xx)$. This implies that for all $k \in \ns$
and $\ell \in \zs$,
$$
[x^\ell ] Y(x)^k= [x^{\ell-k}] Y(\xx)^k= [x^{k-\ell}] Y(x)^k,
\quad  \hbox{that is, } \quad \ct_x (x^{-\ell} Y^k)= \ct_x (x^{\ell-k}
Y^k).
$$
This allows us to rewrite~\eqref{e-Gtfirst} with ``only'' six
terms:
$$
\C(t)=1+ \ct_x \; \left( (\xx^{1}-
x^4)Y + (\xx^{5}-\xx) Y^2 -(\xx^{4}-x^{0}) Y^3\right). 
$$
The
above equation says that $\C(t)$ is the constant term of an
algebraic function. By a very general theory \cite{lipshitz-df},
$\C(t)$ is \emm D-finite,. That is, it satisfies a linear
differential equation with polynomial coefficients. In the next
section, we show that $\C(t)$ satisfies the equation
\eqref{e-dfinite} (or, equivalently, the
P-recurrence~\eqref{e-precursive}). Note that this recurrence can be
easily \emm
guessed, using the {\sc Maple} package {\sc Gfun}: indeed, the
first 15 values of $C_3(n)$ already yield the correct recursion.



\subsection{The Lagrange inversion formula and creative telescoping}
\label{section-lagrange}

From now on, several routes lead to  
the recurrence relation of
Proposition~\ref{main-thm}, depending on how
much software one is willing to use. We present here
the one that we believe to be the shortest. Starting from~\eqref{eq-Y},
the Lagrange inversion formula gives
\beq
\label{a}
[t^n] \ct_x \; ( x^\ell Y^k )= \sum_{j\in \zss} a_n( \ell,k,j)\quad
\hbox{ with }
\quad a_n(\ell,k,j)= \frac k n {n\choose j} {n \choose
  {j+k}} {{2j+k} \choose {j-\ell}}.
\eeq
By convention, the binomial coefficient ${a\choose b}$ is zero
unless $0\le b  \le a$. Hence for $n\ge 1$,
\beq
\label{Ca}
C_3(n)= \sum_{j\in \zss}\big(
a_n(-1,1,j) - a_n(4,1,j)
+ a_n(-5,2,j-1) - a_n(-1,2,j-1)-
a_n(-4,3,j-1) + a_n(0,3,j-1)\big).
\eeq
Of course, we could replace all occurrences of $j-1$ by $j$ in the
above expression, but this results in a bigger final formula.
\begin{Proposition}\label{closed-form}
For $n\ge 1$, the number of $3$-noncrossing partitions of $[n]$ is
$$
C_3(n)=\sum_{j=1}^n
\frac {4 (n-1)!\; n!\; (2j)!}{(j-1)!\;j!\;(j+1)! \;(j+4)! \;(n-j)!\; (n-j+2)!} P(j,n)
$$
with
$$
P(j,n)=
 24 + 18\,n
+ \left( 5-13\,n \right) j
+ \left( 11\,n +20\right) {j}^{2}
+ \left( 10\,n -2 \right) {j}^{3}
+ \left( 4\,n-11 \right) {j}^{4}
-6\,{j}^{5}.
$$
 \end{Proposition}

\medskip
\noindent {{\bf {Proof of the   recurrence relation of
      Proposition~\ref{main-thm}.}} 
We finally apply to the above expression Zeilberger's algorithm for
creative telescoping~\cite[Ch.~6]{AB} (we used the {\sc Maple} package
{\sc Ekhad}). This gives
a recurrence relation for the sequence $C_3(n)$: for $n\ge 2$,
$$
9n ( n-2 ) ( n-
3 )  \left( 4{n}^{2}+15n+17 \right)  C_3 ( n-3 ) - ( n-2 )  \left( 76{n}^{4
}+373{n}^{3}+572{n}^{2}+203n-144 \right) C_3 ( n-2 )
$$
$$+
 ( n+3 )  \left( 44{n}^{4}+189{n}^{3}+227{n}^{2}+30
n-160 \right) C_3 ( n-1 ) = ( n+5 )  ( n+4
 )  ( n+3 )  \left( 4{n}^{2}+7n+6 \right) C_3
 ( n ).
$$
The initial conditions are $C_3(0)=C_3(1)=1$. 
It is then very simple to check that the three term
recursion given in Proposition~\ref{main-thm} satisfies also
the above four term recursion.
This can be rephrased by saying that~\eqref{e-precursive}  is a right
factor of the four term 
recursion obtained via creative telescoping. 
More precisely, applying  to~\eqref{e-precursive} the operator
$$\left( n+1 \right)  \left(
4\,{n}^{2}+39\,n+98 \right) - \left( n+6 \right)  \left(
4\,{n}^{2}+31\,n+63 \right)N ,$$ where $N$ is the shift operator
replacing $n$ by $n+1$, gives the four term
recursion (with $n$ replaced by $n+3$).
 \qed

\medskip
It is worth noting that though the sequence $\{
w_2((1,0),(1,0),2n): n\in \NN  \}$ satisfies a simple recurrence
of order $2$, {\sc Maple} tells us the sequences $\{
q_2((1,0),(1,0),2n): n\in \NN  \}$ and $\{
q_2((1,0),(0,1),2n): n\in \NN  \}$ satisfy more complicated
recurrences of order $3$.


\subsection{Asymptotics}\label{sec:asympt}
Finally, we will  derive from the expression~\eqref{Ca} of $C_3(n)$ the
asymptotic behaviour of this sequence, as stated in Proposition~\ref{main-thm}. 
For any fixed values of $k$ and $\ell$, with $k\ge 1$, we consider the numbers
$a_n(\ell,k,j)$ defined by~\eqref{a}. For the sake of simplicity, we
denote them $a_n(j)$, and introduce the numbers
$$
b_n(j)=  {n\choose j} {n \choose
  {j+k}} {{2j+k} \choose {j-\ell}},
$$
{so that } $ a_n(j)=  k\,  b_n(j) /n.$ Let
$
B_n= \sum_j b_n(j)$. The sum runs over $j\in[\ell_+, n-k]$, where
$\ell_+=\max(\ell,0)$. Below, we often consider $j$ as a \emm real,  variable
running in that interval. The number $b_n(j)$ is then defined in terms
of the Gamma function rather than in terms of
 factorials. 
We will show that $B_n$ admits, for any $N$, an
expansion of the form
\beq\label{Bn-exp}
B_n=9^n \left( \sum_{i= 1}^N \frac{c_i}{n^i}+ O(n^{-N-1})\right),
\eeq
where the coefficients $c_i$ depend on $k$ and $\ell$, and explain how
to obtain these coefficients. We follow the standards steps for
estimating sums of positive terms that are described, for instance,
in~\cite[Section~3]{bender}. 
We begin with a unimodality property of the numbers $b_n(j)$.

\begin{Lemma}[{\bf Unimodality}]\label{lemma-uni}
For  $n$ large enough,  the sequence $b_n(j)$, for $ j \in[\ell_+, n-k]$, is
unimodal, and   its maximum is reached in the interval
$[2n/3-k/2-1/2, 2n/3-k/2-1/3]$. 
\end{Lemma}
\noindent{\bf Proof.} One has
$$
q_n(j):=\frac{b_n(j)}{b_n(j+1)}= \frac 1{(n-j)(n-j-k)}
\frac{(j+1)(j+k+1)}{2j+k+1}
\frac{(j-\ell+1)(j+k+\ell +1)}{2j+k+2}.
$$
Each of these three factors is easily seen to be an increasing function of
$j$ on the relevant interval. Moreover, for $n$ large enough,
$q_n(\ell_+)<1$ while $q_n(n-k-1)>1$.  Let $j_0$
be the smallest value of $j$ such that $q_n(j)\ge 1$. Then $b_n(j)<
b_n(j+1)$ for $j<j_0$ and $b_n(j)\ge b_n(j+1)$ for $j \ge j_0$. We
have thus proved unimodality. Solving  $q_n(x)=1$ for $x$ helps to
locate the mode:
$$
x= \frac{2n} 3  - \frac{5+6k}{12} +O(1/n).
$$
It is then easy to check that $q_n(2n/3-k/2-1/2)<1$ and
$q_n(2n/3-k/2-1/3)>1$ for $n$ large enough.\qed

\smallskip
The second step of the proof reduces the range of summation.
\begin{Lemma}[{\bf A smaller range}]
\label{reduced}
Let $\eps \in (0,1/6)$. Then for all $m$,
$$
B_n=\sum_{|j-2n/3|\le n^{1/2+\eps}}b_n(j)+ o(9^nn^{-m}).
$$  
\end{Lemma}
\noindent{\bf Proof.} Let 
$j=2n/3\pm n^{1/2+\eps}$. The Stirling formula gives
\beq\label{bn-estimate}
b_n(j)= \left(\frac 3 2\right)^{5/2} \frac{9^n}{(\pi n)^{3/2}}
e^{- {9n^{2\eps}/2}}(1+o(1))=o(9^n n^{-m})
\eeq
for all $m$ (the details of the calculation are given in greater
detail below, in the proof of Lemma~\ref{lemma-bn}). Thus by Lemma~\ref{lemma-uni},
$$
\sum_{|j-2n/3|>n^{1/2+\eps}}b_n(j)\le n
\left(b_n(2n/3-n^{1/2+\eps})
+ b_n(2n/3+n^{1/2+\eps})\right)= o(9^nn^{-m}).
$$
The result follows.\qed

\smallskip
A first order  estimate of $b_n(j)$,
 generalizing~\eqref{bn-estimate},  suffices to obtain, upon summing
over $j$, an estimate of $B_n$ of the form~\eqref{Bn-exp} with
$N=1$. However, numerous cancellations occur  when summing the 6 terms
$a_n(\ell,k,j)$ in the expression~\eqref{Ca} of $C_3(n)$. This explains
why we have to work out a longer expansion of the numbers
$b_j(n)$ and $B_n$. 

\begin{Lemma}[{\bf Expansion of $b_j(n)$}]\label{lemma-bn}
  Let $\eps \in (0,1/6)$. Write 
$j=2n/3+ r$ with $r=s\sqrt n$ and
 $|s|\le n^{\eps}$. Then for all $N\ge 1$, 
$$
b_n(j
) = {n\choose j} {n \choose
  {j+k}} {{2j+k} \choose {j-\ell}}=
\left(\frac 3 2\right)^{5/2} \frac{9^n}{(\pi n)^{3/2}}
e^{-{9r^2}/(2n)}
\left(\sum_{i=0}^{N-1}
\frac{c_i(s)}{n^{i/2}}+O(n^{N(3\eps-1/2)})\right)
$$
where  $c_i(s)$ is a  polynomial in $k$, $\ell$ and $s$, of degree
$3i$ in $s$.
Moreover, $c_i$ is an even [odd] function of $s$ if $i$
is even [odd]. In particular, $c_0(s)=1.$
This expansion is uniform in $s$. 
\end{Lemma}
\noindent{\bf Proof.} Note that we simply want to \emm prove the existence, of an
expansion of the above form. The coefficients can be obtained
routinely using (preferably) a computer algebra system. 
In what follows, $(c_i)_{i\ge 0}$ stands for a
sequence of real numbers such that $c_0=1$. The actual value of $c_i$
may change from one formula to another. Similarly, $(c_i(s))_{i\ge 0}$
denotes a sequence of polynomials in $s$ such that $c_0(s)=1$,
having the parity property stated in the lemma.

We start from the Stirling expansion of the Gamma function: for all
$N\ge 1$,
\beq\label{est1}
\Gamma(n+1)=n^n \sqrt{2\pi n}\ e^{-n}
\left(\sum_{i=0}^{N-1}\frac{c_i}{n^i} + O(n^{-N})\right).
\eeq
This gives, for $j=2n/3+r$, with $r=s\sqrt n$, and for any $N$,
\beq\label{est2}
\Gamma(j+1)=2 j^j  \sqrt{\frac{\pi n}3} e^{-j}
\left(\sum_{i=0}^{N-1}\frac{c_i(s)}{n^{i/2}} + O(n^{N(\eps-1/2)})\right) 
\eeq
for some polynomials $c_i(s)$ of degree $i$ in $s$. This estimate is
uniform in $s$.  Similarly,
\beq\label{est3}
\Gamma(n-j+1)= (n-j)^{n-j} \sqrt{\frac{2\pi n} 3}\ e^{j-n}
\left(\sum_{i=0}^{N-1}\frac{c_i(s)}{n^{i/2}} + O(n^{N(\eps-1/2)})\right),
\eeq
for some polynomials $c_i(s)$ of degree $i$ in $s$.
Now
\begin{eqnarray*}
\log \frac {n^n}{j^j (n-j)^{n-j}}&= &\log \frac{3^n}{2^j} -\frac
     {9r^2}{4n}-\sum_{i\ge 2} \frac{3^i\, r^{i+1}}{i(i+1) n^i}
     \left(1-(-1/2)^i\right)\\
&= &\log \frac{3^n}{2^j} -\frac
     {9r^2}{4n}-\sum_{i= 1}^{N-1} \frac{c_i(s)}{n^{i/2}} +
     O(n^{2\eps+N(\eps-1/2)}),
\end{eqnarray*}
for some polynomials $c_i(s)$ of degree $i+2$ in $s$. Observe that
$(i+2)/(i/2)\le 6$ for $i\ge 1$.
Hence 
\begin{eqnarray}
\frac {n^n}{j^j (n-j)^{n-j}}&= & \frac{3^n}{2^j} e^{-{9r^2}/({4n})}
\exp\left( -\sum_{i= 1}^{N-1} \frac{c_i(s)}{n^{i/2}} +
     O(n^{2\eps+N(\eps-1/2)})\right)\nonumber\\
&=& \frac{3^n}{2^j}
e^{-{9r^2}/(4n)}\left(\sum_{i=0}^{N-1}\frac{c_i(s)}{n^{i/2}} +
O(n^{N(3\eps-1/2)})\right),\label{est4}
\end{eqnarray}
for polynomials $c_i(s)$ of degree $3i$ in $s$.
Putting together (\ref{est1}--\ref{est4}), one obtains, uniformly in
$s$:
\beq\label{est5}
{n\choose j}= \frac 3{2\sqrt {\pi n}}  \frac{3^n}{2^j} e^{-{9r^2}/(4n)}\left(\sum_{i=0}^{N-1}\frac{c_i(s)}{n^{i/2}} +
O(n^{N(3\eps-1/2)})\right),
\eeq
for polynomials $c_i(s)$ of degree $3i$ in $s$.

Similarly,
\beq\label{est6}
{n\choose {j+k}}= \frac 3{2\sqrt {\pi n}}  \frac{3^n}{2^{j+k}} 
e^{-{9r^2}/({4n})}
\left(\sum_{i=0}^{N-1}\frac{c_i(s)}{n^{i/2}} +
O(n^{N(3\eps-1/2)})\right),
\eeq
with the same degree condition on the polynomials $c_i(s)$. Finally,
\beq\label{est7}
{{2j+k}\choose {j-l}}= \sqrt{\frac{3}{2 \pi n}} 2^{2j+k} \left(\sum_{i=0}^{N-1}\frac{c_i(s)}{n^{i/2}} +
O(n^{N(\eps-1/2)})\right)
\eeq
for polynomials $c_i(s)$ of degree $i$ in $s$.
Putting together  (\ref{est5}--\ref{est7}), we obtain the estimate of
$b_n(j)$ given in the lemma.
\qed

\smallskip
It remains to sum our estimates of $b_n(j)$ for values of $j$ such that
$|j-2n/3|\le n^{1/2+\eps}$.
\begin{Proposition}[{\bf Expansion of $B_n$}]\label{propo-Bn} For all $N\ge 1$,
$$
B_n= \sum_j  {n\choose j} {n \choose
  {j+k}} {{2j+k} \choose {j-\ell}}
= 
\left(\frac 3 {2}\right)^{5/2} \frac{9^n}{\pi ^{3/2}n}
\left(\sum_{i=0}^{N-1}\frac{1}{n^i}
\int_{\rs} 
c_{2i}(s) e^{- {9s^2}/{2}}
ds
+O(n^{-N})
\right)
$$
where the polynomials $c_i(s)$, depending on $k$ and $\ell$,  are
those of Lemma~{\em\ref{lemma-bn}.} 
In particular, as $c_0(s)=1$,
$$
B_n(j)=  \frac{3^{3/2}}{4\pi} \frac{9^n}{n} \left( 1 +O(1/n)\right).
$$
\end{Proposition}
\noindent{\bf Proof.} We start from  Lemma~\ref{reduced}, and combine it with
the uniform expansion of
Lemma~\ref{lemma-bn}. We need to estimate sums of the following type, for  $i\in \ns$:
$$
\sum_{|j-2n/3|\le n^{1/2+\eps}} \left(\frac{j-2n/3}{\sqrt
    n}\right)^i e^{-9(j-2n/3)^2/(2n)}.
$$
Using the Euler-MacLaurin summation
formula~\cite[Eq.~(5.62)]{odlyzko-handbook}, one 
obtains, for all $m\ge 1$,
$$
\frac 1 {\sqrt n}\sum_{|j-2n/3|\le n^{1/2+\eps}} \left(\frac{j-2n/3}{\sqrt
    n}\right)^i e^{-9(j-2n/3)^2/(2n)}= \int_{\rs} s^i e^{-9s^2/2} ds +
  o(n^{-m}).
$$
The above integral vanishes if $i$ is odd, 
and can be expressed in
terms of the Gamma function otherwise. This gives the estimate of
Proposition~\ref{propo-Bn}, but with a rest of order
$n^{N(6\eps-1)}$. However, this expansion is valid for all $N$, and
for all $\eps>0$. From this, the rest can be seen to be of order
$n^{-N}$.
\qed

\smallskip
With the strategy described above (and {\sc Maple}...), we have obtained
the expansion of $B_n$  to the order  $n^{-6}$. Given that
$a_n(j)\equiv a_n(\ell,k,j)=k\, b_n(j)/n$, this gives the expansion of $\sum_j
a_n(\ell,k,j)$ to the order  $n^{-7}$:
$$
\sum_j a_n(\ell,k,j)=   \frac{3^{3/2}}{4\pi} \frac{9^n}{n^2} \left(
\sum _{i=0}^5 \frac {c_i}{(4n)^i} + O(1/n^6)\right).
$$
The coefficients of this expansion are too big to be reported here
for generic values of $k$ and $\ell$ (apart from $c_0=k$). We simply give the values of
$c_i$ for the 6 pairs $(\ell,k)$ that are involved in the
expression~\eqref{Ca} of $C_3(n)$:

$$
\begin{array}{c|cccccc|ccccccc}
(\ell,k) & c_0& c_1&c_2&c_3&c_4 & c_5
\\
\hline
(-1,1) & 1& -7& 37& -184& 871& -4087 & + \\
(4,1) &  1& -127& 8317& -381904& 14034391& -444126847 &-  \\
(-5,2) & 2& -230& 13682& -573416& 19338062& -564941270 &+ \\
(-1,2) &  2& -38& 434& -4136& 36302& -305558 &-\\
(-4,3) & 3& -237& 9831& -293664& 7227813& -157405197&-\\
 (0,3) & 3& -165& 4863& -106104& 1959573& -32693205&+
\end{array}
$$
Each pair $(\ell,k)$ contributes in the expression of $C_3(n)$ with a
weight $\pm1$, depending on the sign  indicated on the corresponding
line of the above table. One observes that the first 5 terms cancel,
which leaves
$$
C_3(n)=   \frac{3^{3/2}}{4\pi} \frac{9^n}{n^2} \left(
\frac {4199040}{(4n)^5} + O(1/n^6)\right).
$$
This completes the proof of Proposition~\ref{main-thm}.\qed
 
\subsection{Other starting points}

So far, we have focused on the enumeration of $Q_2$-vacillating walks
starting at $(1,0)$. However, our approach works just as well for
other starting points, or for
combinations of starting points. Let $A(x,y)$ be the generating
function of starting points. For instance, $A(x,y)=x^iy^j$
corresponds to starting at $(i,j)$, while $A(x,y)=1/(1-x)$ corresponds
to starting anywhere on the $x$-axis.

Then a calculation similar to that of Section~\ref{sec-functional} gives the
following functional equation
$$
K(x,y;t)F^A(x,y;t) = (x+y+xy)A(x,y) - xH^A(x;t)-yV^A(y;t), 
$$
where $F^A,$ $H^A,$ and $V^A$ are the analogues of $F$, $H$ and $V$.
Specializing $A(x,y)$ to $x$ gives back \eqref{e-Fkernel}.
The kernel method of Section~\ref{sec-kernel} now gives
\begin{multline*}
xH^A(x)+\xx V^A(\xx)=\\
\left( x+Y+xY \right) A \left( x,Y \right)
-
\left(\xx Y+ Y+\xx Y^2 \right) A \left(\xx Y, Y \right)
+
\left( \xx Y +\xx +\xx ^2 Y\right) A \left(\xx Y, \xx \right)
.
\end{multline*}
Thus if $A$ is a rational function (and in particular if $A(x,y)=x^iy^j$),
then $H^A(x)$, $V^A(x)$, and $F^A(x,y)$ are all $D$-finite.
%

\section{Partitions with no enhanced 3-crossing}
\label{sec-hesitant}

Our approach for counting 3-noncrossing partitions can be easily
adapted to the enumeration of partitions avoiding \emm enhanced,
3-crossings. As discussed
in~\cite{crossings-nestings},  partitions of $[n]$ avoiding
enhanced 
$k+1$-crossings are in bijection with \emm hesitating, tableaux of
height at most $k$. In turn, these hesitating tableaux are in
one-to-one correspondence with certain $W_k$-\emph{hesitating}
lattice walks. A hesitating lattice walk satisfies the following
walking rules: when pairing every two steps from the beginning,
each pair of steps has one of the following three types: i) a stay
step followed by an $e_i$ step, ii) a $-e_i$ step followed by a
stay step, iii) an $e_i$ step followed by a $-e_j$ step.

Then partitions of $[n]$ avoiding enhanced $k+1$-crossings are in
bijection with $W_k$-hesitating walks of length $2n$ starting and ending at
$(k-1, \ldots, 2,1, 0)$. As before, $W_k$ denotes the Weyl chamber
$\{\, (a_1,a_2,\dots ,a_{k})\in \ZZ^{k}: a_1> a_2> \cdots >
a_{k}\ge 0 \,\}$.

For convenience, all notations we used for vacillating walks will be recycled
for hesitating  walks. Thus
$\delta=(k-1,k-2,\dots,0)$, and for two lattice points $\lambda$
and $\mu$ in $W_k$, we denote by $w_k(\lambda,\mu,n)$ (resp.
$q_k(\lambda,\mu,n)$) the number of $W_k$-hesitating (resp.
$Q_k$-hesitating) lattice walks of length $n$ starting at
$\lambda$ and ending at $\mu$. A careful investigation
shows that the reflection principle ``works'' for hesitating lattice
walks.

\begin{Proposition}
\label{propo-reflection-h} For any starting and ending points
$\lambda$ and
  $\mu$ in $W_k$, the number of $W_k$-hesitating walks going from
  $\lambda$ to $\mu$ can be expressed in terms of the number of
  $Q_k$-hesitating walks as follows:
$$
w_k(\lambda,\mu,n)=\sum_{\pi\in \sy_k} (-1)^\pi
q_k(\lambda,\pi(\mu),n),
$$
 where $(-1)^\pi$ is the sign of $\pi$ and
$\pi(\mu_1,\mu_2,\dots,\mu_k)=(\mu_{\pi(1)},\mu_{\pi(2)},\dots,
\mu_{\pi(k)})$.
\end{Proposition}

\begin{proof}
The key property here is that the set of pairs of steps that are
allowed for a hesitating walk is left invariant when reflecting either
one, or both steps with respect to the hyperplane $x_i=x_j$. This is
clearly seen by abbreviating the three types of step pairs as $(0, +),
(-,0)$ and $(+,-)$, and recalling that the above reflection exchanges
the $i$th and $j$th coordinates (equivalently, the unit vectors $e_i$
and $e_j$). The rest of the argument copies the proof of  Proposition
\ref{propo-reflection}.
\end{proof}

Let us now focus on the case $k=2$. The connection between
$W_2$-hesitating walks and partitions avoiding enhanced 3-crossings
entails
\beq
\label{cross-hesit}
E_3(n)= w_2((1,0),(1,0),2n)= q_2((1,0),(1,0),2n)-q_2((1,0),(0,1),2n).
\eeq
Let us  write
a functional equation counting $Q_2$-hesitating walks that start from
$(1,0)$.
Let $a_{i,j}(n):=q_2((1,0),(i,j),2n)$ be the number of
such walks having length $2n$ and ending at $(i,j)$. Let
$$
F(x,y;t)=\sum_{i,j,n}a_{i,j}(2n)x^iy^jt^{n}
$$
 be the associated generating
function. Then by appending to an (even length) walk an allowed pair
of steps, we obtain the following functional equation:
\begin{multline*}
\left(x+y+\xx+\yy+(x+y)(\xx+\yy)\right)tF(x,y;t)\\ =
F(x,y;t) -x + H(x;t)(\yy+x\yy )t+V(y;t)(\xx+\xx y)t,
\end{multline*}
where $H(x;t)$ (resp.~$V(y;t)$) is the generating function of even
lattice walks ending on the $x$-axis (resp.~$y$-axis).
This functional equation can be rewritten as
\begin{align}\label{e-H-F}
K(x,y;t) F(x,y;t)= x^2y-x (1+x )tH(x;t)-y(1+y)tV(y;t),
\end{align}
where $K(x,y;t)$ is the kernel given by
$$
K(x,y;t)=xy-t(1+x)(1+y)(x+y).
$$
From now on, we will very often omit the variable $t$ in our
notation.

Let us now solve~\eqref{e-H-F}.
 First fix $x$, and consider the kernel as a quadratic polynomial
in $y$. Only one of its roots, denoted $Y$ below, is a formal
series in $t$:
$$Y= \frac
{1-t\xx (1+x)^2 - \sqrt {\left(1-t\xx(1+x)^2\right)^2-4t\xx(1+x)^2} }
{ 2\left( 1+\xx \right) t}=O(t). 
$$
The coefficients of this series are Laurent polynomials in $x$, as
is easily seen from the equation 
\beq \label{eq-H-YF}
Y=t(1+\xx)(1+Y)(x+Y).
\eeq
The second root of the kernel is $Y_1=x/Y= O(t^{-1})$, and  the
expression $F(x,Y_1)$  is not well-defined.

Observe that the new kernel only differs from~\eqref{kernel} by a term
$txy$. Hence, as in the case of vacillating walks, the product of the
roots is $x$,  and the kernel is symmetric in $x$
and $y$. This implies that the diagram of the roots, obtained by
taking iteratively conjugates, is still given by
Figure~\ref{diagram}. Again, the 3  pairs of power series that are
framed can be legally substituted for $(x,y)$ in the
functional equation~\Ref{e-H-F}. We thus obtain:
\begin{eqnarray}
x(1+x)t H(x)+ Y(1+Y)t V(Y)&=& x^2Y,
\label{e-H-eq1}\\
\xx Y(1+\xx Y)t H(\xx Y)+ Y(1+Y)t V(Y)&=&\xx^2Y^3,
\label{e-H-eq2}\\
\xx Y(1+\xx Y)t H(\xx Y)+ \xx(1+\xx)t V(\xx)&=&\xx^3Y^2.
\label{e-H-eq3}
\end{eqnarray}
Now \eqref{e-H-eq1}$-$\eqref{e-H-eq2}$+$\eqref{e-H-eq3} gives
\begin{align*}
x(1+x)t H \left( x \right) + \xx(1+\xx)t V \left( \xx\right)
=x^2Y-\xx^2Y^3+\xx^3Y^2.
\end{align*}
By~\eqref{eq-H-YF}, the series $Y/t/(1+x)$ is a formal series in $t$
with coefficients in $\qs[x,\xx]$. Thus we can divide the above
identity by $t(1+x)$, and then extract the positive and negative
parts.
\begin{Proposition}
\label{propo-H-walks} The series $H(x)$ and $V(y)$, which  count
$Q_2$-hesitating walks of even length ending on the $x$-axis and
on the $y$-axis, satisfy
\begin{align*}
x H(x)= \pt_x\; \frac Y{t(1+x)} (x^2-\xx^2Y^2+\xx^3Y),\\
\xx^2 V(\xx)= \nt_x\; \frac Y{t(1+x)} (x^2-\xx^2Y^2+\xx^3Y).
\end{align*}
\end{Proposition}
\noindent
Let us now return to the number $E_3(n)$ of
partitions of $[n]$ avoiding enhanced $3$-crossings, given
by~\eqref{cross-hesit}. The generating
function $\mathcal{E}(t)$ of these numbers is:
\begin{align}\label{e-H-Gtfirst}
\E(t)=\ct_x \; \frac {Y(\xx^2-x^3)}{t(1+x)} (x^2-\xx^2Y^2+\xx^3Y).
\end{align}
Observe that, again, $Y(x)=x Y(\xx)$. Therefore, for all $k\in \ns $
and $\ell \in \zs$,
$$
\ct_x\left( \frac{\xx^\ell Y^k}{t(1+x)}\right) =
\ct_x\left( \frac{x^{\ell-k+1} Y^k}{t(1+x)}\right).
$$
This allows us to  rewrite~\eqref{e-H-Gtfirst} with only non-negative
powers of $x$:
\beq\label{e-Gtsecond-H}
\E(t)= \ct_x \; \frac Y{t(1+x)} \left( 1-x^5 -(x^2-x)Y^2 +(x^4-1)Y
\right).
\eeq
The above equation shows that $\E(t)$ is
the constant term of an algebraic function. It is thus \emm
D-finite,. Let us now compute  a linear differential equation it
satisfies (equivalently, a P-recursion for its coefficients).

Starting from~\eqref{eq-H-YF}, the Lagrange inversion formula
gives
$$
[t^n] \ct_x \;  \frac {x^\ell Y^k}{t(1+x)}= \sum_{j\in \zss} a_n(
\ell,k,j)\quad \hbox{ with }
\quad a_n(\ell,k,j)= \frac k {n+1}
{{n+1}\choose j} {{n+1} \choose
  {j+k}} {{n} \choose {j-\ell}}.
$$
From~\eqref{e-Gtsecond-H}, we obtain
\beq\label{E-expr}
E_3(n)= \sum_{j\in\zss}\left(
a_n(0,1,j)-a_n(5,1,j)-a_n(2,3,j)+a_n(1,3,j)+a_n(4,2,j)-a_n(0,2,j)\right).
\eeq
This gives an explicit (but not so simple) expression of $E_3(n)$, to
which we apply  Zeilberger's algorithm for creative
telescoping. This proves that, for $n\ge 1$,
\begin{multline*}
8n \left( n-1 \right)  \left( n-2 \right) 
E_3 ( n-3 )
+3 \left( 5{n}^{2}+17n+8 \right)  \left( n-1 \right) E_3 (
n-2)  \\
+3 \left( n+1 \right)  \left( 2n+5
 \right)  \left( n+4 \right) E_3 ( n-1 )
- \left( n+6
 \right)  \left( n+5 \right)  \left( n+4 \right) E_3 ( n )
 =0,
\end{multline*}
with initial condition $E_3(0)=1$. It is then straightforward to
check that the sequence defined  in Proposition \ref{main-thm-H}
also satisfies the above P-recursion. More precisely, applying the
operator $(n+2)+(n+7)N$ to the recursion of Proposition
\ref{main-thm-H} gives the above four term recursion.

The study of the aymptotic behaviour of $E_3(n)$  parallels what we did for
$C_3(n)$. The maximum of $a_n(\ell,k,j)$ is now reached for $j\sim
n/2$ rather than $2n/3$. Using the same notations as in
Section~\ref{sec:asympt}, we obtain
$$
B_n(j) \sim \frac {8^{n+1}}{\sqrt 3 \pi n },$$
but again, numerous cancellations occur when we sum the 6 required
estimates, so as to obtain the estimate of Proposition \ref{main-thm-H}.

\section{Final comments}
It is natural to ask whether for any $k$, the sequence $C_k(n)$ that
counts $k$-noncrossing partitions of $[n]$ is P-recursive. Our opinion
is that this is unlikely, at least for $k=4$. This is based on the
following observations:
\begin{enumerate}
  \item We have written a functional equation for $Q_3$-vacillating
  walks, with kernel
$$K(x,y,z;t)=
1-t(1+x+y+z)(1+\xx+\yy+\zz).$$
Using this equation, we have computed the first 100 numbers in the
sequence $C_4(n)$. This is sufficient for the {\sc Maple} package {\sc
  Gfun} to discover a P-recursion of order 8 with coefficients of
degree 8, if it exists. But no such recursion has been found.
\item
Let us solve the above kernel in $z$. The two roots  $Z_0$ and $Z_1$ are
related by
\beq
\label{Z-relation}
Z_0 Z_1= \frac{1+x+y}{1+\xx+\yy}.
\eeq
Since the kernel is symmetric in $x$, $y$ and $z$, the diagram of the
roots, obtained by taking conjugates, is generated by the
transformations  $\Phi_i$ for $i=1,2,3$, where
$$\Phi_3(x,y,z)=\left(x,y,\zz \, \frac{1+x+y}{1+\xx+\yy}\right)$$
 and $\Phi_1$ and
$\Phi_2$ are defined similarly. But these transformations now generate an \emm
infinite, diagram. There exist in the literature a few signs 
indicating that a finite
(resp. infinite) diagram is related to a D-finite (resp. non-D-finite)
\gf. First, a number of equations with a finite diagram have been
solved, and shown to have a D-finite solution~\cite{bousquet-versailles,bousquet-motifs,bousquet-kreweras,mishna}. Then, the only
problem with an infinite diagram that has been thoroughly studied has
been proved to be
non-D-finite~\cite{bousquet-petkovsek-knight}.
Finally,  the conjectural
link between finite diagrams and D-finite series is confirmed by a systematic
numerical study of walks in the quarter plane~\cite{mishna}.
\end{enumerate}

The above paragraphs can be copied verbatim for $W_3$-hesitating walks
and partitions avoiding enhanced $4$-crossings. The kernel is now
$$
K(x,y,z)=
1-t(x+y+z+\xx+\yy+\zz+(x+y+z)(\xx+\yy+\zz)),
$$
but the roots $Z_0$ and $Z_1$ are still related by~\eqref{Z-relation}.

\smallskip
As recalled in the introduction, the sequence $M_k(n)$ that counts
$k$-noncrossing \emm matchings, of $[n]$ (that is, partitions in which
all blocks have size 2) is D-finite for all $k$. More precisely, the
associated exponential \gf,
$$
\M_k(t)=\sum_n M_k(n) \frac {t^{2n}}{(2n)!}
$$
is given by~\cite{grabiner}:
$$
\M_k(t)=\det \left( I_{i-j}(2t)-I_{i+j}(2t) \right)_{1\le i,j \le k-1},$$
where
$$
I_n(2t)=\sum_{j\ge 0} \frac{t^{n+2j}}{j!(n+j)!}$$
is the hyperbolic Bessel function of the first kind of order
$n$. 
The existence of such a closed form implies that $\M_k(t)$ is
D-finite~\cite{lipshitz-df}. The specialization to matchings of
the bijection between partitions and vacillating
walks results in a bijection between $k+1$-noncrossing
matchings and \emm oscillating tableaux, of height at most $k$,
or, equivalently, $W_k$-oscillating walks. These walks can take
any (positive or negative) unit step $\pm e_i$, without any parity
restriction. The kernel of the equation ruling the enumeration of
such walks is simply
$$
K(x_1, \ldots , x_k)=1-t( x_1+ \cdots + x_k +\xx_1+ \cdots +
\xx_k).
$$
The diagram of the roots, generated by the $\Phi_i$, for $1\le i
\le k$, where $\Phi_i(x_1, \ldots , x_k)=(x_1, \ldots, x_{i-1},
\xx_i,x_{i+1}, \ldots , x_k),$ is now finite (the group of
transformations $\Phi_i$ being itself isomorphic to
$(\zs/2\zs)^k$).  This, again, confirms the possible connection
between finite diagrams and D-finite series.

\bigskip \noindent
{\bf Acknowledgements:} We are grateful to Anders
       Bj\"{o}rner and Richard Stanley for inviting us to the
``Algebraic Combinatorics'' program at the Institut Mittag-Leffler
in Spring 2005, during which  this work was done. We also thank
Christian Krattenthaler for explaining us how the asymptotics of our
sequences could be worked out completely.  MBM  was partially
supported by the European Commission's IHRP 
  Programme, grant HPRN-CT-2001-00272, ``Algebraic Combinatorics in
  Europe''. GX would like to thank Richard Stanley and Rosena Du for
introducing him to this topic.

\bibliographystyle{plain}
\bibliography{biblio}

\end{document}